\documentclass[12pt]{article}
\input tcilatex
\usepackage{amsfonts,epsf}

\newcommand{\psdiag}[3]{\hspace{1mm}\raisebox{-#1mm}{\epsfysize#2mm
\epsffile{#3.eps}}\hspace{1mm}}
\newcommand{\psbild}[2]{\vspace{3mm}\\ \mbox{\epsfysize=#1mm\epsffile{#2.eps}}
\vspace{3mm}}

\begin{document}

\author{Rui Pedro Carpentier\\
\\
{\small\it Departamento de Matem\'{a}tica and Centro de
Matem\'{a}tica Aplicada}\\
{\small\it  Instituto Superior T\'{e}cnico}\\
{\small\it Avenida Rovisco Pais, 1049-001 Lisboa}\\
{\small\it Portugal}}
\title{From planar graphs to embedded graphs - a new approach to Kauffman and
Vogel's polynomial}
\date{3rd May, 2000}
\maketitle

\begin{abstract}

In \cite{4} Kauffman and Vogel constructed a rigid vertex regular isotopy
invariant for unoriented
four-valent graphs embedded in three dimensional space. It assigns to each
embedded graph $G$ a
polynomial, denoted $[G]$,  in three variables, $A$, $B$ and $a$, satisfying
the skein relations:
$$
[\psdiag{2}{6}{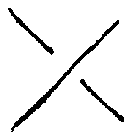}]=A [\psdiag{2}{6}{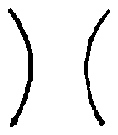}]+B
[\psdiag{2}{6}{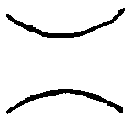}]+
[\psdiag{2}{6}{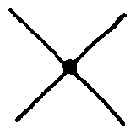}]$$
$$
[\psdiag{2}{6}{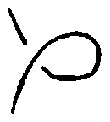}]= a[\psdiag{2}{6}{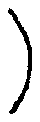}],    \hspace{1cm}
[\psdiag{2}{6}{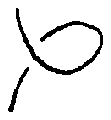}]= a^{-1} [\psdiag{2}{6}{straight}]
$$
and is defined in terms of a state-sum and the Dubrovnik polynomial for
links.
Using the graphical calculus of \cite{4} it is shown that the
polynomial of a planar graph can
be calculated recursively from that of planar graphs with less vertices,
which also allows the
polynomial of an embedded graph to be calculated without resorting to links.
The same approach is
used to give a direct proof
of uniqueness of the (normalized) polynomial restricted to
planar graphs. In the case
$B=A^{-1}$ and $a=A$, it is proved that for a planar graph  $G$ we have
$[G]=2^{c-1}(-A-A^{-1})^v $, where $c$ is the number of connected components
of $G$ and $v$ is the number of vertices of $G$. As a corollary, a
necessary, but not sufficient,
condition is obtained  for an embedded graph to be ambient isotopic to a
planar graph. In an appendix
it is shown that, given a polynomial for planar graphs satisfying the
graphical calculus, and imposing
the first skein relation above, the polynomial extends to a rigid vertex
regular isotopy invariant for
embedded graphs, satisfying the remaining skein relations. Thus, when
existence of the planar
polynomial is guaranteed, this provides a direct way, not depending on
results for the Dubrovnik
polynomial, to show consistency of the polynomial for embedded graphs.

\end{abstract}

\section{Introduction.}

An {\it ambient isotopy} for 4-valent rigid vertex embedded graphs may be
regarded \cite{1} as a finite sequence of the following moves:

$$
\begin{array}{c}
I) \psdiag{2}{6}{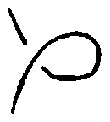} \rightleftharpoons \psdiag{2}{6}{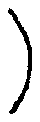} \rightleftharpoons \psdiag{2}{6}{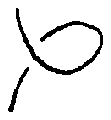}\\ 
II) \psdiag{2}{6}{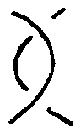} \rightleftharpoons \psdiag{2}{6}{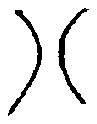}\\
III) \psdiag{2}{6}{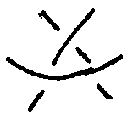} \rightleftharpoons \psdiag{2}{6}{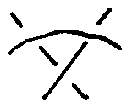}\ \ ,
\psdiag{2}{6}{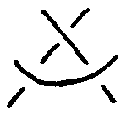} \rightleftharpoons \psdiag{2}{6}{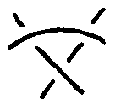}\\
IV) \psdiag{2}{6}{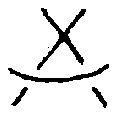} \rightleftharpoons \psdiag{2}{6}{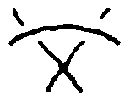}\ \ ,\psdiag{2}{6}{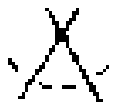} \rightleftharpoons\psdiag{2}{6}{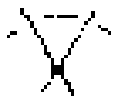}\\
V)\psdiag{2}{6}{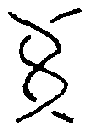} \rightleftharpoons \psdiag{2}{6}{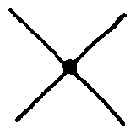} \rightleftharpoons \psdiag{2}{6}{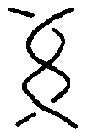}\end{array}
$$

If the first move is not used in the sequence it is called a {\it regular
isotopy}.

The {\it Dubrovnik polynomial }\cite{2,3} is the knot polynomial in 2-variables 
$a,z$ invariant for regular isotopies that satisfies the axioms:

$$
\begin{array}{c}
i)D_{\psdiag{2}{4}{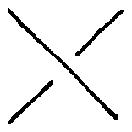}}-D_{\psdiag{2}{4}{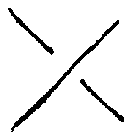}}=z(D_{\psdiag{2}{4}{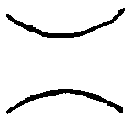}}-D_{\psdiag{2}{4}{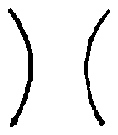}})
\\ 
ii)D_{\psdiag{2}{4}{f1}}=aD_{\psdiag{2}{4}{g1}},D_{\psdiag{2}{4}{h1}}=a^{-1}D_{\psdiag{2}{4}{g1}}
\\ 
iii)D_{\psdiag{2}{4}{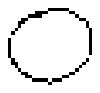}}=1 
\end{array}
$$

Kauffman and Vogel \cite{4} construct a 3-variable polynomial for 4-valent
rigid vertex embedded graphs which is invariant under regular isotopies, by using
the skein relation $[\psdiag{2}{6}{a1} ]=A[\psdiag{2}{6}{d1}
]+B[\psdiag{2}{6}{e1} ]+[\psdiag{2}{6}{c1}
]$ and the Dubrovnik polynomial with $z=A-B$.

It is constructed as follows. Let $G$ be a 4-valent rigid vertex embedded graph
diagram.

(1) Choose a marker \psdiag{2}{6}{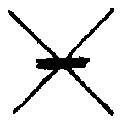}\ \ or \psdiag{2}{6}{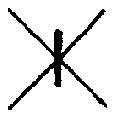}\ \ at each vertex
of the graph.

(2) Let ${\cal L}$ be the set of links obtained by replacing:

$$
\begin{array}{c}
\psdiag{2}{6}{j2} \longmapsto \psdiag{2}{6}{b1} or\psdiag{2}{6}{e1} or  \psdiag{2}{6}{d1}\\
\psdiag{2}{6}{i2} \longmapsto \psdiag{2}{6}{a1} or\psdiag{2}{6}{d1} or \psdiag{2}{6}{e1}
\end{array}
$$

(3) Define $[G]=\stackunder{L\in {\cal L}}{\sum }(-A)^{i(L)}(-B)^{j(L)}D_L$
where $i(L)$ is the number of replacements of type \psdiag{2}{6}{j2}\ \
$\longmapsto $
\psdiag{2}{6}{e1}\ \ , $j(L)$ is the number of replacements of type
\psdiag{2}{6}{j2}\ \ $%
\longmapsto $\psdiag{2}{6}{d1}\ \ and $D_L$ is the Dubrovnik polynomial of $L$
with $%
z=A-B$.

Kauffman and Vogel \cite{4} prove that this polynomial is well-defined and
invariant under regular isotopy. It generalizes both the bracket
polynomial (corresponding to the case $B=A^{-1}$ and $a=-A^3$) and the Yamada
polynomial 
\cite{5} (corresponding to the case $B=A^{-1}$ and $a=A^2$).

\section{Graphical Calculus.}

In \cite{4} the following graphical calculus is proved.

\begin{theorem}
For 4-valent graph diagrams, differing only in the marked local picture, we
have the following identities:

\begin{description}
\item  $[\psdiag{2}{6}{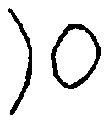} ]=\mu [\psdiag{2}{6}{g1}]$

\item  $[\psdiag{2}{6}{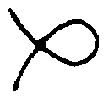} ]={\cal O}[\psdiag{2}{6}{g1} ]$

\item  $[\psdiag{2}{6}{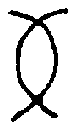} ]=(1-AB)[\psdiag{2}{6}{k1} ]+\gamma
[\psdiag{2}{6}{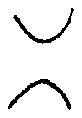} ]-(A+B)[\psdiag{2}{6}{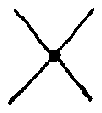} ]$

\item  $[\psdiag{2}{6}{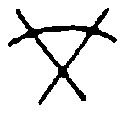}]-[\psdiag{2}{6}{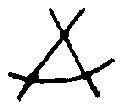}
]=AB([\psdiag{2}{6}{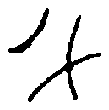}]-[\psdiag{2}{6}{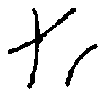}]+[\psdiag{2}{6}{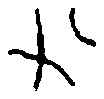}]-[\psdiag{2}{6}{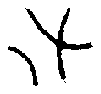}
]+[\psdiag{2}{6}{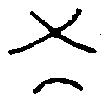}]-[\psdiag{2}{6}{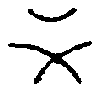}])+\xi
([\psdiag{2}{6}{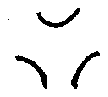}]-[\psdiag{2}{6}{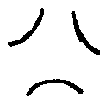}])$
\end{description}

where
$$
\begin{array}{c}
\mu = 
\frac{a-a^{-1}}{A-B}+1 \\ {\cal O}=\frac{Aa^{-1}-Ba}{A-B}-(A+B) \\ \gamma = 
\frac{B^2a-A^2a^{-1}}{A-B}+AB \\ \xi =\frac{B^3a-A^3a^{-1}}{A-B} 
\end{array}
$$
\end{theorem}

These identities allow us to calculate the polynomial of a planar
4-valent graph from the  polynomials of other planar graphs with less vertices,
if the graph contains one of the local configurations \psdiag{2}{6}{l2} or
\psdiag{2}{6}{m2}
 or if it contains one of these configurations after a sequence of moves of
the type \psdiag{2}{6}{a3} $%
\longrightarrow $\psdiag{2}{6}{b3} . This is in fact always the case, as the
following
result shows.

\begin{lemma}
If a connected planar 4- valent graph has at least one vertex and does not
contain either of the local
configurations \psdiag{2}{6}{l2} or \psdiag{2}{6}{m2} then it is possible,
by a
sequence of moves of the type \psdiag{2}{6}{a3} $\longrightarrow
$\psdiag{2}{6}{b3} , to
turn it into a new planar 4-valent graph containing the local configuration
\psdiag{2}{6}{m2}
.
\end{lemma}

\TeXButton{Proof}{\proof}First we prove that the graph contains a global
configuration of the type:%

\centerline{\psbild{25}{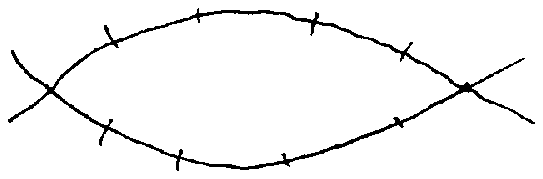} \hspace{2.5cm} \mbox{}}

Let $v$ be a vertex of the graph. We choose an edge that departs from $v$,
and construct a walk
which starts along this edge and carries on straight ahead at each vertex it
meets. In other words, when 
the walk gets to a vertex it continues along the opposite edge, as opposed to
turning to the left or right.

\centerline{\psbild{25}{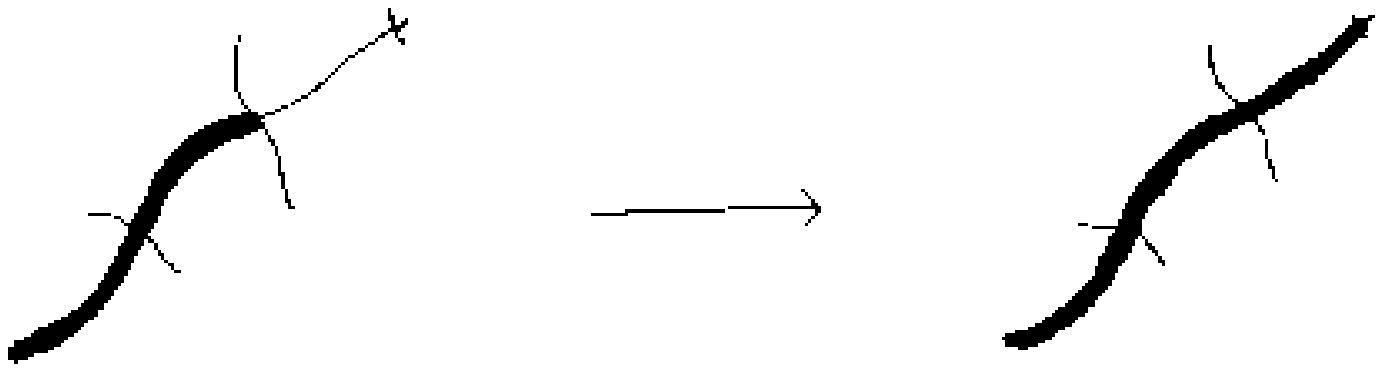}}

At a certain stage the walk returns to a vertex it has already visited producing
one of the following
three global configurations:

$$
\begin{array}{c}
\psdiag{5}{20}{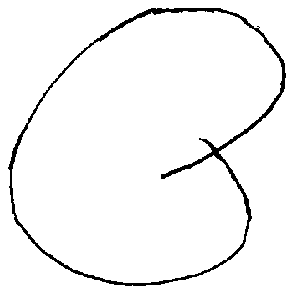} \hspace{0.5cm} {\rm or} \hspace{0.5cm} 
\psdiag{5}{20}{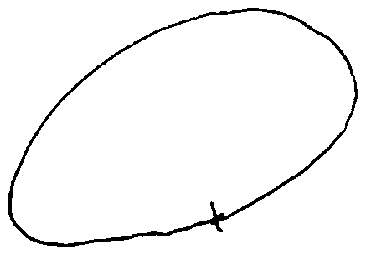} \hspace{0.5cm} {\rm or} \hspace{0.5cm} \psdiag{5}{20}{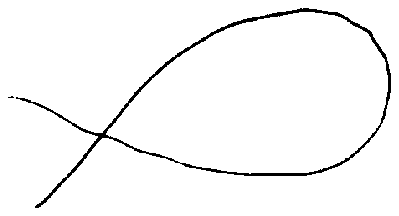}
\end{array}
$$

In the first case we prolong the walk so that we 
either obtain the third global configuration above, if the walk terminates or
crosses itself before leaving the enclosed region:
$$
\psdiag{5}{20}{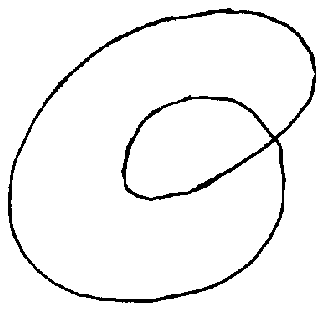} \hspace{1cm} {\rm or} \hspace{1cm} \psdiag{5}{20}{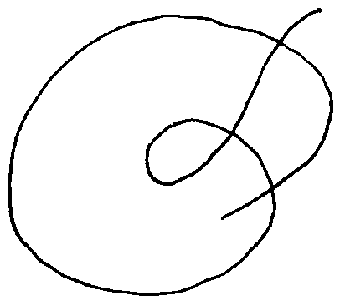}
$$

or the desired configuration,
if the walk leaves the enclosed region without crossing itself before it does so:

\centerline{\psbild{25}{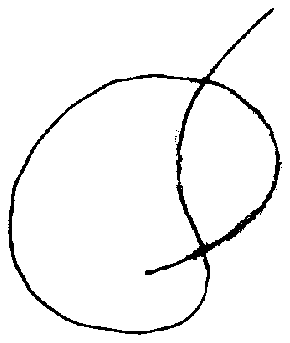}\hspace{2.5cm} \mbox{}}

In the second case, starting at vertex $v$,  we follow the transversal walk to
the first walk into the enclosed region, and in this way we either obtain the
third global configuration above:%

\centerline{\psbild{25}{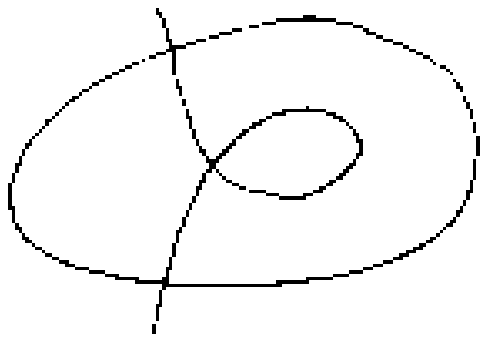}\hspace{2.5cm} \mbox{}}

or the desired configuration: 

\centerline{\psbild{25}{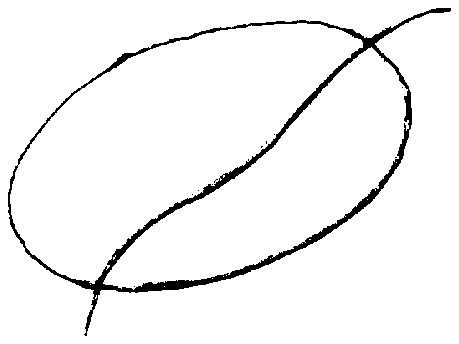}\hspace{2.5cm} \mbox{}}

In the third case, we take the set of all configurations of this type in the
graph and choose one that is minimal (in the sense that it does not
enclose another configuration of the same type).
Since, by hypothesis, the graph does not contain loops, we can choose any
walk that crosses this minimal configuration, and in this way get the desired
configuration:%

\centerline{\psbild{25}{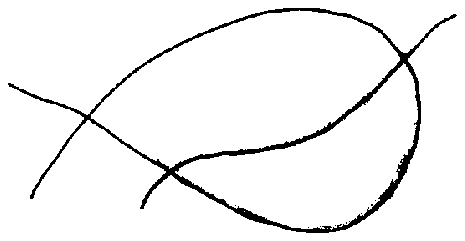}\hspace{2.5cm} \mbox{}}

Now, since the graph contains configurations of the type 

\centerline{\psbild{25}{lem7}\hspace{2.5cm} \mbox{}}

we can choose one that is minimal (in the sense that it does not contain
another configuration of same type inside).

This configuration does not contain configurations of the second or third type
 inside it,
since it is minimal. Thus there are no walks starting inside the configuration 
which 
do not eventually cross to the outside, and any walk that goes in through one of
the sides of
the configuration goes out through the other side without intersecting
itself.

Suppose, first, that the configuration has vertices inside it.  Start a walk at
one of these
vertices and extend it until it crosses to the outside. At the last vertex before
this walk crosses 
to the outside, take the transversal walk that leaves the configuration through
the same side.

\centerline{\psbild{25}{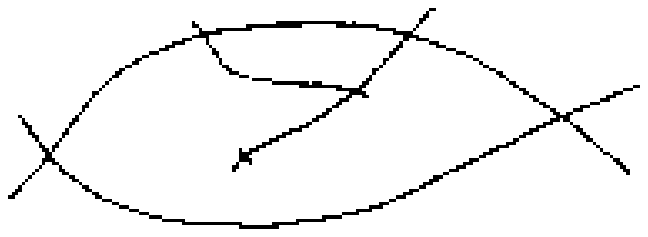}\hspace{2.5cm} \mbox{}}

Thus, we get a global triangle adjacent to one side of the configuration.

At the last vertex of the second walk before it crosses to the outside, take the
transversal walk
that goes into the triangle. This walk must leave  the triangle
through the side of the configuration, because if the walk went out through
the same side of the triangle as it went in,
the configuration would not be minimal, and the third side has no vertices
available for it to cross.
$$
\psdiag{2}{20}{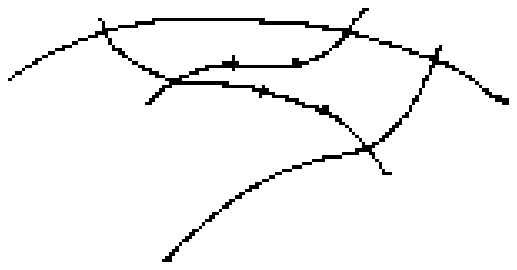} \hspace{1cm} {\rm ,} \hspace{1cm} \psdiag{2}{20}{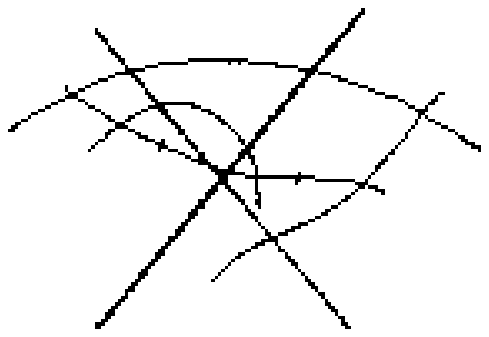}
$$

We iterate this last step so as to obtain a finite chain of
triangles $T_1\supset T_2\supset ...\supset T_n$, all adjacent to one
side of the configuration.%

\centerline{\psbild{25}{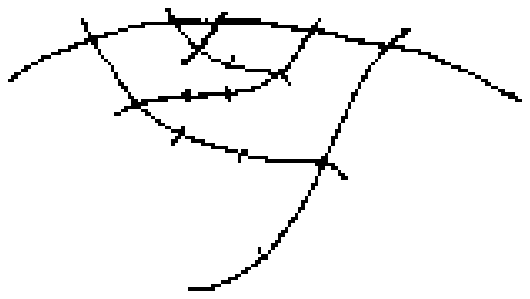}\hspace{2.5cm} \mbox{}}

The last of these is a simple triangle.Thus, we can remove one vertex from inside
the
configuration by a move of the type \psdiag{2}{6}{a3} $\rightleftharpoons $\psdiag{2}{6}{b3} :%

\centerline{\psbild{25}{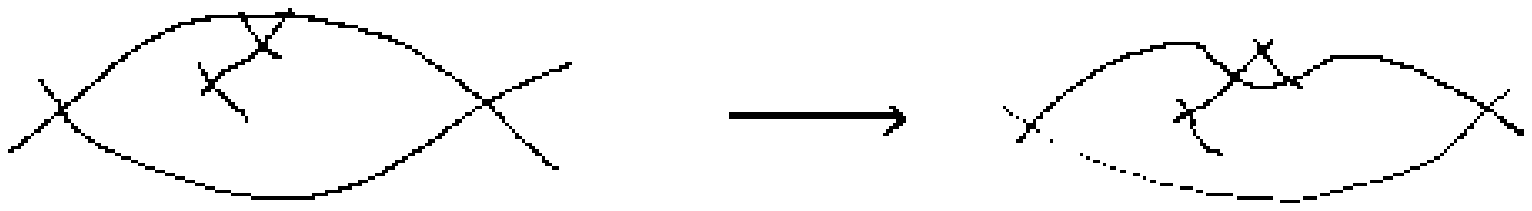}}

Applying this method recursively we can remove all vertices from inside the
configuration, which takes on the following form:%

\centerline{\psbild{25}{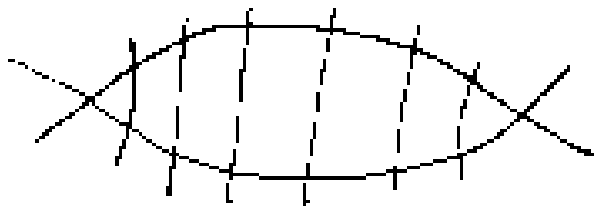}\hspace{2.5cm} \mbox{}}

Then we can apply the move \psdiag{2}{6}{a3} $\rightleftharpoons $\psdiag{2}{6}{b3}
successively removing all edges inside the configuration. At the end, we
obtain the desired configuration \psdiag{2}{6}{m2} (with no edges or vertices
inside).%

\centerline{\psbild{25}{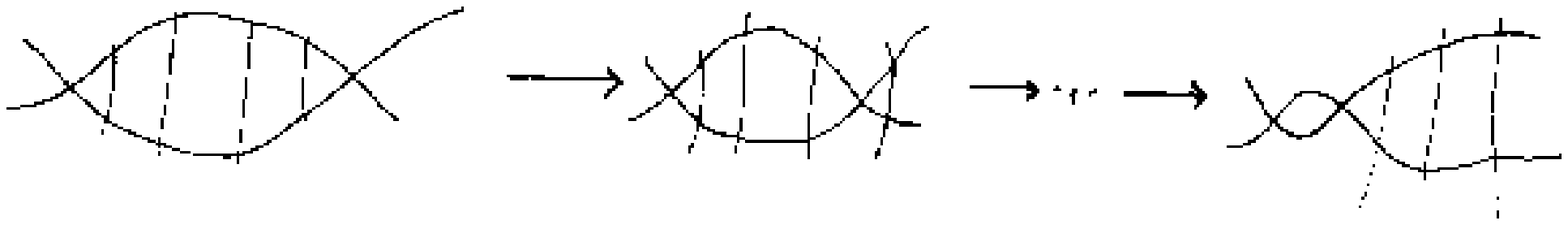}} 

\TeXButton{End Proof}{\endproof}

With this graphical calculus we can calculate the polynomial of planar
graphs without resorting to links.
We note that, since the polynomial is invariant under regular isotopies, the
position of the various connected components of the graph does not
matter: e.g. $[\psdiag{2}{6}{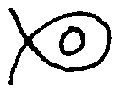} ]=[\psdiag{2}{6}{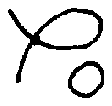} ]$ .
In fact, we have the following result:

\begin{theorem}
There is a unique polynomial for planar 4-valent graphs that satisfies
the graphical calculus and takes value 1 for the unknot.
\end{theorem}

\TeXButton{Proof}{\proof}Let $[.]$ and $[.]^{*}$ be two polynomials satisfying
the
conditions of the theorem. We will prove, by induction on the number of
vertices, that $[G]=[G]^{*}$ for any planar 4-valent graph $G$.

If $G$ has no vertices then $G$ is an disjoint union of unknots, and hence, by
the first identity of graphical calculus, $[G]=[G]^{*}=\mu ^{c-1}$ where $c$
is the number of unknots in $G$.

By the induction hypothesis, we suppose that $[.]$ and $[.]^{*}$ are equal for
graphs with up to $n$ vertices.

Let $G$ be a graph with $n+1$ vertices.

If $G$ contains a local configuration of the type \psdiag{2}{6}{l2} (or
\psdiag{2}{6}{m2} )
then, by the second (or third) identity of the graphical calculus and the
induction
hypothesis, we get $[G]=[G]^{*}$.

Otherwise, by lemma 2 it is possible to produce a sequence of graphs $%
G\longrightarrow G_1\longrightarrow ....\longrightarrow G_m$ using moves of the
type \psdiag{2}{6}{a3}\ \ $\rightleftharpoons $\psdiag{2}{6}{b3}\,  where $G_m$
contains a
configuration of the type \psdiag{2}{6}{m2}\ . Then, using the last identity of
the 
graphical calculus and the induction hypothesis, we conclude that $%
[G]-[G_1]=[G]^{*}-[G_1]^{*}$, $[G_1]-[G_2]=[G_1]^{*}-[G_2]^{*}$, ...,$%
[G_{m-1}]-[G_m]=[G_{m-1}]^{*}-[G_m]^{*}$, and hence
$[G]-[G_m]=[G]^{*}-[G_m]^{*}$. 
Now $[G_m]=[G_m]^{*}$ because $G_m$ contains a configuration of the type 
\psdiag{2}{6}{m2}, and thus we conclude that $[G]=[G]^{*}$.\TeXButton{End
Proof}{\endproof}

Moreover, we can calculate the polynomial of a 4-valent rigid vertex
embedded graph using this calculus in the following way:

Let $G$ be a 4-valent rigid vertex embedded graph. Let ${\cal P}$ be the
set of planar graphs obtained by replacing 
$$
\psdiag{2}{6}{a1} \longmapsto \psdiag{2}{6}{d1} ,\psdiag{2}{6}{e1} or
\psdiag{2}{6}{c1}
$$

Then $[G]=\stackunder{P\in {\cal P}}{\sum }A^{i(P)}B^{j(P)}[P]$ where $i(P)$
is the number of replacements of type \psdiag{2}{6}{a1} $\longmapsto
$\psdiag{2}{6}{d1}
and $j(P)$ is the number of replacements of type \psdiag{2}{6}{a1} $\longmapsto $%
\psdiag{2}{6}{e1} .

It is also possible to show that the invariance of the polynomial under regular
isotopies is a
consequence of the identities of theorem 1 (see appendix).

\section{The case when $B=A^{-1}$ and $a=A$.}

If we make the choice $a=A$ and $B=A^{-1}$ then: 
$$
\begin{array}{c}
\mu =2, \\ 
{\cal O}=-A-A^{-1}, \\ \gamma =0, \\ 
\xi =-A-A^{-1}. 
\end{array}
$$

In this case the graphical calculus takes the form:%
$$
\begin{array}{c}
\lbrack 
\psdiag{2}{6}{k2} ]=2[\psdiag{2}{6}{g1} ], \\ \lbrack 
\psdiag{2}{6}{l2} ]=-(A+A^{-1})[\psdiag{2}{6}{g1} ], \\ \lbrack 
\psdiag{2}{6}{m2} ]=-(A+A^{-1})[\psdiag{2}{6}{l9} ], \\ \lbrack
\psdiag{2}{6}{b3}]-[\psdiag{2}{6}{a3}]=[%
\psdiag{2}{6}{c3}]-[\psdiag{2}{6}{e3}]+[\psdiag{2}{6}{g3}]-[\psdiag{2}{6}{f3}]+
[\psdiag{2}{6}{h3}]-[\psdiag{2}{6}{i3}]\\
-(A+A^{-1})([\psdiag{2}{6}{j3}]-[\psdiag{2}{6}{l4}]). 
\end{array}
$$

\begin{theorem}
In the case $a=A$ and $B=A^{-1}$ for any 4-valent planar graph $G$ we have $%
[G]=2^{c-1}(-A-A^{-1})^v$, where $c$ is the number of connected components
of $G$ and $v$ is the number of vertices of $G.$
\end{theorem}

\TeXButton{Proof}{\proof}By Theorem 3 we only have to prove that $[\psdiag{2}{6}{m6}
]=1$ (that is obvious) and that the polynomial $2^{c-1}(-A-A^{-1})^v$
satisfies the graphical calculus with $B=A^{-1}$ and $a=A$.

The three first identities are evidently satisfied.

Now, let us look at the last identity. This identity is
invariant under the 6-dihedral group $D_6$. If we remove from the graphs in
the equation the configurations in which they differ and label the free
ends $1$ to $6$ we have, up to
the symmetry of the equation, the following distinct cases:

\begin{description}
\item[1.]  $1,2,3,4,5,6$ are in the same connected component $\psdiag{2}{10}{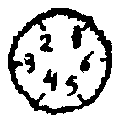}$;

\item[2.]    $1,2$ are in one connected component and $3,4,5,6$ in another
$\psdiag{2}{10}{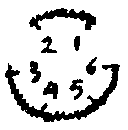}$;

\item[3.]    The free ends group into connected components as:  $ \{
1,2\},\{3,4\},\{5,6\}$ $\psdiag{2}{10}{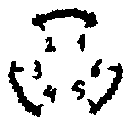}$;

\item[4.]   The free ends group into connected components as: 
$\{1,2\},\{3,6\},\{4,5\}$ $\psdiag{2}{10}{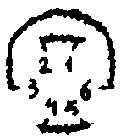}$.
\end{description}

In case 1 we have $[\psdiag{2}{6}{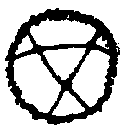}]=[\psdiag{2}{6}{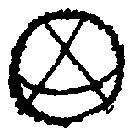}]$,
$[\psdiag{2}{6}{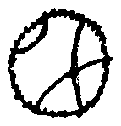}]=[\psdiag{2}{6}{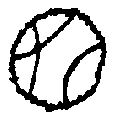}]$, $[%
\psdiag{2}{6}{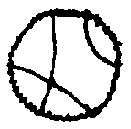}]=[\psdiag{2}{6}{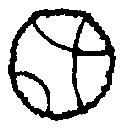}]$,
$[\psdiag{2}{6}{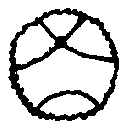}]=[\psdiag{2}{6}{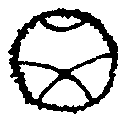}]$ and
$[\psdiag{2}{6}{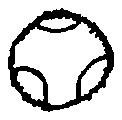}]=[\psdiag{2}{6}{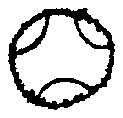}]$,
and thus the identity is satisfied.

In case 2 we have $[\psdiag{2}{6}{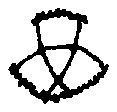}]=[\psdiag{2}{6}{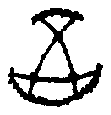}]$,
$[\psdiag{2}{6}{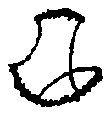}]=[%
\psdiag{2}{6}{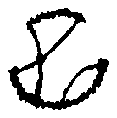}]$, $[\psdiag{2}{6}{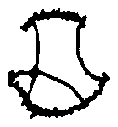}]=[\psdiag{2}{6}{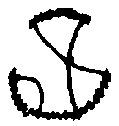}]$,
$2[\psdiag{2}{6}{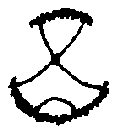}]=[\psdiag{2}{6}{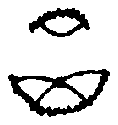}]$ and $[\psdiag{2}{6}{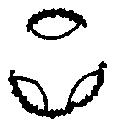}]=2[\psdiag{2}{6}{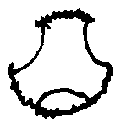}]$,
 and thus the identity
is $0=[\psdiag{2}{6}{h8}]-(-A-A^{-1})[\psdiag{2}{6}{k8}]=(-A-A^{-1})[\psdiag{2}{6}{k8}]
-(-A-A^{-1})[\psdiag{2}{6}{k8}]=0$, i.e. the identity is satisfied.

In case 3 we have $[\psdiag{2}{6}{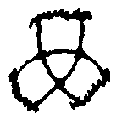}]=[\psdiag{2}{6}{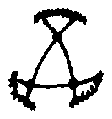}]$,
$2[\psdiag{2}{6}{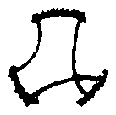}]=[%
\psdiag{2}{6}{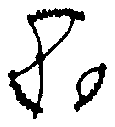}]$, $2[\psdiag{2}{6}{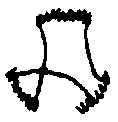}]=[\psdiag{2}{6}{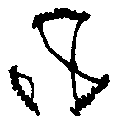}]$,
$2[\psdiag{2}{6}{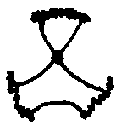}]=[\psdiag{2}{6}{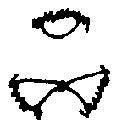}]$ and $[\psdiag{2}{6}{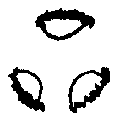}]=4[\psdiag{2}{6}{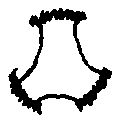}]$,
 and thus the identity is $0=[%
\psdiag{2}{6}{c9}]+[\psdiag{2}{6}{g9}]+[\psdiag{2}{6}{h9}]-3(-A-A^{-1})[\psdiag{2}{6}{k9}]=
3(-A-A^{-1})[\psdiag{2}{6}{k9}]-3(-A-A^{-1})[\psdiag{2}{6}{k9}]=0$, i.e.
the identity is
satisfied.

In case 4 we have $[\psdiag{2}{6}{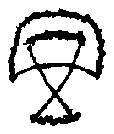}]=[\psdiag{2}{6}{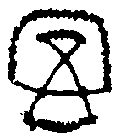}]$,
$[\psdiag{2}{6}{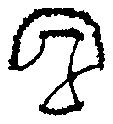}]=[%
\psdiag{2}{6}{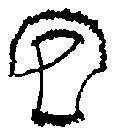}]$, $[\psdiag{2}{6}{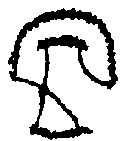}]=[\psdiag{2}{6}{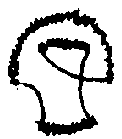}]$,
$[\psdiag{2}{6}{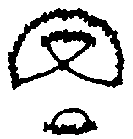}]=[\psdiag{2}{6}{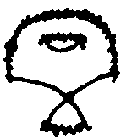}]$ and $[\psdiag{2}{6}{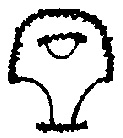}]=[\psdiag{2}{6}{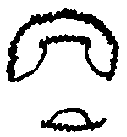}]$,
 and thus the identity
is satisfied.%
\TeXButton{End Proof}{\endproof}

 If we replace the vertices of a graph $G$ by crossings we obtain a link $L$ .
A {\it knot-theoretic circuit} of $G$ is 
a closed walk on $G$ which corresponds to 
a link
component of $L$. We define the {\it writhe} of a knot-theoretic circuit $D$
as $w(D)=\stackunder{c\in C(D)}{\sum }\epsilon (c)$ where $C(D)$ is the set
of  self-crossings of $D$ and $\epsilon (c)$ is the sign $ (\pm 1)$ of the
crossing 
$c$ for a chosen orientation (and independent of this choice). The {\it twisting
number} 
$t(G)$ of $G$ is the
sum of the writhes over all knot-theoretic circuits of $G$.

\begin{corollary}
A necessary condition for a 4-valent rigid vertex embedded graph $G$ be isotopic
to a planar graph is that
the polynomial of $G$ with $B=A^{-1}$ and $a=A$ is $%
[G]=2^{c-1}(-A-A^{-1})^vA^{t(G)}$ where $c$ is the number of connected
components, $v$ is the number of vertices and $t(G)$ is the twisting number
of $G$.
\end{corollary}

This corollary is a consequence of the fact that the polynomial $%
a^{-t(G)}[G] $ is an isotopy invariant.

Unfortunately, this condition is not sufficient as we see in the following
example.

$$
G:\psdiag{14}{28}{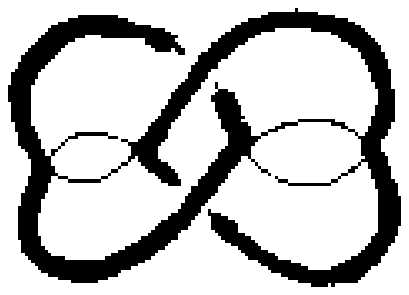} \hspace{1cm} c=1,t(G)=2,v=4.
$$

\begin{example}
The graph in the picture is not planar, since it contains two disjoint linked
cycles (shown bold in the picture), 
 but its polynomial is $%
A^2(-A-A^{-1})^4=2^{c-1}A^{t(G)}(-A-A^{-1})^v$.
\end{example}

\begin{corollary}
If a graph diagram $G$ has only one crossing, and the removal of this crossing
does not
change the number of connected components, then the polynomial of $G$
vanishes in the case $B=A^{-1}$ and $a=A$, and thus $G$ is not isotopic to a
planar graph.
\end{corollary}

\TeXButton{Proof}{\proof}$[G]=A[G_1]+A^{-1}[G_2]+[G_3]$ where $G_1$ is the
graph obtained by the change \psdiag{2}{6}{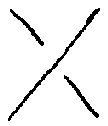} $\longrightarrow
$\psdiag{2}{6}{k1} , $%
G_2$ is the graph obtained by the change \psdiag{2}{6}{k11} $\longrightarrow$
\psdiag{2}{6}{d10} and $G_3$ is the graph obtained by the change
\psdiag{2}{6}{k11} $%
\longrightarrow$ \psdiag{2}{6}{l9}.
Since the number of connected components of $G$ does not change when the crossing

is removed,  $G_1$,$G_2$
and $G_3$ have the same number of connected components, and therefore
$[G_1]=[G_2]$
and $[G_3]=(-A-A^{-1})[G_1]$. Thus $[G]=(A+A^{-1})[G_1]+(-A-A^{-1})[G_1]=0$.%
\TeXButton{End Proof}{\endproof}

\section{Comments.}

We can prove that the graphical calculus implies the invariance of the
polynomial for regular isotopies (see the appendix). Thus, if we can prove
the consistence of the graphical calculus without using the invariance of the
polynomial, then we don't need to use the Dubrovnik polynomial to prove the
invariance of the polynomial. This is achieved in the case $B=A^{-1}$ and $%
a=A$ because Theorem 4 gives a polynomial that is consistent with the graphical
calculus, and which is unique by Theorem 3.

We can also see that this choice of variables ($B=A^{-1}$ and $a=A$) is the
most interesting case, assuming that for planar graphs the
polynomial is of the form $[G]=p^{c-1}q^v$ , where $p$ and $q$ are two
polynomials, $c$ is the number of connected components and $v$ is the
number of vertices of $G$. In fact, by the two first identities of the
graphical calculus we get $p=\mu $ and $q={\cal O}$. By a reasoning
analogous to the proof of Theorem 4, the third identity gives us the
following equations:%
$$
\begin{array}{c}
{\cal O}^2=1-AB+\gamma -(A+B){\cal O} \\ {\cal O}^2=(1-AB)\mu +\gamma -(A+B)%
{\cal O} \\ {\cal O}^2=1-AB+\gamma \mu -(A+B){\cal O} 
\end{array}
$$

If $\mu \neq 1$ then we have $AB=1$ and $\gamma =0$, thus $a=A$ or $a=-A^3$.
The second case implies ${\cal O}=0$, and thus we get a generalization of the
bracket polynomial that vanishes for graphs with at least one vertex. If $\mu
=1$ it is easy to prove, by induction on the number of crossings, that for any
4-valent rigid vertex embedded graph $G$ we have $[G]=(A+B+{\cal O})^{cr}{\cal
O}^v$,
where $cr$ is the number of the crossings of the embedding and $v$ is the
number of vertices of $G$, and we also have $A+B+{\cal O}=a$ and $%
a=a^{-1}$. Thus, in this case, the invariant $a^{-t(G)}[G]$ is almost trivial: 
  $a^{-t(G)}[G]={\cal O}^v$ if $t(G)$
and $cr$ have the same parity (both even or both odd) and $a^{-t(G)}[G]=a{\cal
O}^v$ otherwise.

{\bf Acknowledgment} - I wish to  thank Prof. Roger Picken who encouraged me to write 
this paper and helped me to improve it.

\appendix
\section{Appendix.}

We show that the graphical calculus implies the regular isotopy invariance
of the polynomial.

\begin{proposition} 
If $[G]$ is a polynomial for 4-valent graph diagrams that satisfies the
graphical calculus, and the skein relation $[\psdiag{2}{6}{a1}
]=A[\psdiag{2}{6}{d1}]+B[\psdiag{2}{6}{e1}]+[\psdiag{2}{6}{c1}]$ holds,
 then it is invariant
under regular isotopies for 4-valent rigid vertex embedded graphs  and
satisfies the identities $[\psdiag{2}{6}{f1}]=a[\psdiag{2}{6}{g1}\ ]$ and
$[\psdiag{2}{6}{h1}]=a^{-1}[\psdiag{2}{6}{g1}]$.
\end{proposition}

\TeXButton{Proof}{\proof}1. $[\psdiag{2}{6}{f1}]=
A[\psdiag{2}{6}{k2}]+B[\psdiag{2}{6}{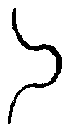}]+[\psdiag{2}{6}{l2}]\\
=A\mu [\psdiag{2}{6}{g1}]+B[\psdiag{2}{6}{g1}]+{\cal O}[\psdiag{2}{6}{g1}]\\
=a[\psdiag{2}{6}{g1}]\\ $ 
 and $[\psdiag{2}{6}{h1} ]=A[\psdiag{2}{6}{m11}]+B[\psdiag{2}{6}{k2}]+[\psdiag{2}{6}{l2}]\\
=A[\psdiag{2}{6}{g1}]+B\mu [\psdiag{2}{6}{g1}]+{\cal O}[\psdiag{2}{6}{g1}]\\
=a^{-1}[\psdiag{2}{6}{g1}]\\ $
 , since ${\cal O}+A\mu +B=a$ and ${\cal O}+A+B\mu =a^{-1}$.

2. To prove invariance under regular isotopy, it is enough to show the invariance of the 
polynomial under the generalized Reidemeister moves $II)- V)$.

$II) \quad [\psdiag{2}{6}{j1}]=A[\psdiag{2}{6}{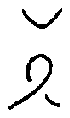}]+B[\psdiag{2}{6}{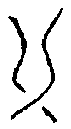}]
+[\psdiag{2}{6}{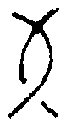}]\\
=Aa^{-1}[\psdiag{2}{6}{d10}]+B[\psdiag{2}{6}{a1}]+[\psdiag{2}{6}{d12}]\\
=Aa^{-1}[\psdiag{2}{6}{d10}]+BA[\psdiag{2}{6}{k1}]
+B^2[\psdiag{2}{6}{d10}]+B[\psdiag{2}{6}{l9}
]+A[\psdiag{2}{6}{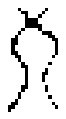}]+B[\psdiag{2}{6}{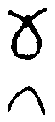}]+[\psdiag{2}{6}{m2}]\\
=AB[\psdiag{2}{6}{k1}]+(Aa^{-1}+B^2+B{\cal
O})[\psdiag{2}{6}{d10}]+(A+B)[\psdiag{2}{6}{l9}]+[\psdiag{2}{6}{m2}]\\
=AB[\psdiag{2}{6}{k1}]-\gamma [\psdiag{2}{6}{d10}
]+(A+B)[\psdiag{2}{6}{l9}]+[%
\psdiag{2}{6}{m2}]=[\psdiag{2}{6}{k1}]$ , by the third identity of the
graphical calculus.

$IV) \quad [\psdiag{2}{6}{f2} ]=A[\psdiag{2}{6}{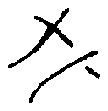}]+B[\psdiag{2}{6}{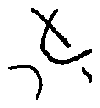}
]+[\psdiag{2}{6}{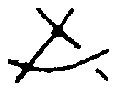} ]=A^2[\psdiag{2}{6}{e3} ]+AB[%
\psdiag{2}{6}{h3}]+A[\psdiag{2}{6}{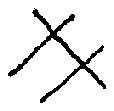} ]+BA[\psdiag{2}{6}{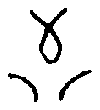}
]+B^2[\psdiag{2}{6}{f3}]+B[\psdiag{2}{6}{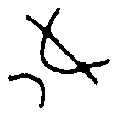} ]+A[%
\psdiag{2}{6}{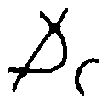}]+B[\psdiag{2}{6}{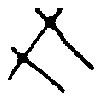}]+[\psdiag{2}{6}{a3}]\\
=A^2[\psdiag{2}{6}{e3}]+AB[\psdiag{2}{6}{h3}]+A[%
\psdiag{2}{6}{m10}]+BA{\cal O}[\psdiag{2}{6}{j3}]+B^2[\psdiag{2}{6}{f3}
]+B(1-AB)[\psdiag{2}{6}{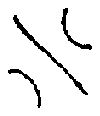}
]+B\gamma [\psdiag{2}{6}{j3}]-B(A+B)[\psdiag{2}{6}{f3}
]+A(1-AB)[\psdiag{2}{6}{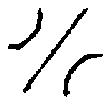}]+A\gamma [%
\psdiag{2}{6}{j3}]-A(A+B)[\psdiag{2}{6}{e3}] +B[\psdiag{2}{6}{b11}]+[\psdiag{2}{6}{a3}]\\
=-AB([\psdiag{2}{6}{e3}]-[%
\psdiag{2}{6}{h3}]+[\psdiag{2}{6}{f3}])+(BA{\cal O}+A\gamma +B\gamma
)[\psdiag{2}{6}{j3}]+B(1-AB)[%
\psdiag{2}{6}{i12}]+A(1-AB)[\psdiag{2}{6}{k12}]+A[\psdiag{2}{6}{m10}
]+B[\psdiag{2}{6}{b11}]+[\psdiag{2}{6}{a3}]\\
=-AB([\psdiag{2}{6}{e3}]-[\psdiag{2}{6}{h3}]+[\psdiag{2}{6}{f3}])+\xi
[\psdiag{2}{6}{j3}]+B(1-AB)[%
\psdiag{2}{6}{i12}]+A(1-AB)[\psdiag{2}{6}{k12}]+A[\psdiag{2}{6}{m10}
]+B[\psdiag{2}{6}{b11}]+[\psdiag{2}{6}{a3}].$

Rotating the equation by 180$^o$ we obtain:

$[\psdiag{2}{6}{g2}]=-AB([\psdiag{2}{6}{c3}]-[\psdiag{2}{6}{i3}
]+[\psdiag{2}{6}{g3}])+\xi [\psdiag{2}{6}{l4}
]+B(1-AB)[\psdiag{2}{6}{i12}]+A(1-AB)[\psdiag{2}{6}{k12}
]+A[\psdiag{2}{6}{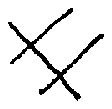}]+B[\psdiag{2}{6}{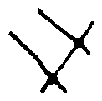}]+[%
\psdiag{2}{6}{b3}]$. Note that the skein relation is invariant under rotations.

Thus $[\psdiag{2}{6}{g2}]-[\psdiag{2}{6}{f2}]=-AB([\psdiag{2}{6}{c3}
]-[\psdiag{2}{6}{e3} ]+[\psdiag{2}{6}{h3}]-[%
\psdiag{2}{6}{i3}]+[\psdiag{2}{6}{g3}]-[\psdiag{2}{6}{f3}])-\xi
([\psdiag{2}{6}{j3}]-[\psdiag{2}{6}{l4}
])+B(1-AB)([\psdiag{2}{6}{i12}]-[\psdiag{2}{6}{i12}
])+A(1-AB)([\psdiag{2}{6}{k12}]-[\psdiag{2}{6}{k12}])+A([%
\psdiag{2}{6}{m10r}]-[\psdiag{2}{6}{m10}])+B([\psdiag{2}{6}{b11r}
]-[\psdiag{2}{6}{b11}])+[\psdiag{2}{6}{b3}]-[\psdiag{2}{6}{a3}]=0$.

$III) \quad [\psdiag{2}{6}{m1}]=A[\psdiag{2}{6}{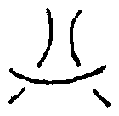}]+B[\psdiag{2}{6}{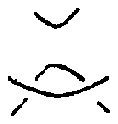}
]+[\psdiag{2}{6}{f2}]\\
=A[\psdiag{2}{6}{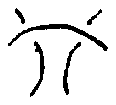}]+B[\psdiag{2}{6}{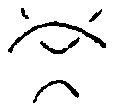}]+[\psdiag{2}{6}{g2}]
=[\psdiag{2}{6}{a2}]$.

$V) \quad [\psdiag{2}{6}{e2}]=-A[\psdiag{2}{6}{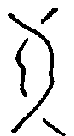}]-B[\psdiag{2}{6}{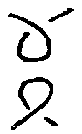}
]+[\psdiag{2}{6}{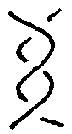}]\\
=-A[\psdiag{2}{6}{k1}]-B[%
\psdiag{2}{6}{d10}]+[\psdiag{2}{6}{k11}]=[\psdiag{2}{6}{l9}
]$.\TeXButton{End Proof}{\endproof}

Since the polynomial of an embedded graph can be given as a weighted sum 
of  polynomials of planar graphs, we only need to assume that the graphical
calculus holds for planar graphs.

\end{document}